\documentclass[smallextended]{svjour3}       
\smartqed  
\usepackage{graphicx}
\usepackage{epstopdf}
\usepackage{amsmath}
\usepackage{amssymb}
\usepackage{color,url}
\usepackage{algorithm,algorithmic}
\usepackage{subfigure}
\usepackage{multirow}
\usepackage{rotating}
\usepackage{tabularx}
\usepackage{booktabs}
\usepackage{lscape}

\usepackage{hyperref}
\hypersetup{hypertex=true,
	colorlinks=true,
	linkcolor=blue,
	anchorcolor=blue,
	citecolor=blue}

\allowdisplaybreaks[4]
\newcommand{\tr}{^{\sf T}}

\begin{document}

\title{Analytic analysis of the worst-case complexity of the gradient method with exact line search and the Polyak stepsize
}

\titlerunning{Analytic analysis of the worst-case complexity of the gradient method ...}        

\authorrunning{Ya-Kui Huang, Hou-Duo Qi} 

\author{Ya-Kui Huang        \and
	Hou-Duo Qi 
}


\institute{
	Ya-Kui Huang 
			\at
	Institute of Mathematics, Hebei University of Technology, Tianjin 300401, China\\
	\email{hyk@hebut.edu.cn}
	\and
	Hou-Duo Qi \at
	Department of Data Science and Artificial Intelligence and 
	Department of Applied Mathematics, The Hong Kong Polytechnic University, Hung Hom, Hong Kong\\
	\email{houduo.qi@polyu.edu.hk}
}

\date{Received: date / Accepted: date}

\maketitle

\begin{abstract}
We give a novel analytic analysis of the worst-case complexity of the gradient method with exact line search and the Polyak stepsize, respectively, which previously could only be established  by computer-assisted proof. Our analysis is based on studying the linear convergence of a family of gradient methods, whose stepsizes include the one determined by exact line search and the Polyak stepsize as special instances.   The asymptotic behavior of the considered family is also investigated which shows that the gradient method with the Polyak stepsize will zigzag in a two-dimensional subspace spanned by the two eigenvectors corresponding to the largest and smallest eigenvalues of the Hessian.
\keywords{Gradient methods \and Exact line search \and Polyak stepsize \and Worst-case complexity}
\subclass{90C25 \and 90C25 \and 90C30}
\end{abstract}

\section{Introduction}
\label{intro}
We focus on the gradient method for solving unconstrained optimization problem
\begin{equation}\label{equncpro}
	\min_{x\in\mathbb{R}^n}~f(x),
\end{equation}
where $f: \mathbb{R}^n \to \mathbb{R}$ is a real-valued, convex, and continuously differentiable function. As is known, the gradient method updates iterates by
\begin{equation}\label{eqitr}
	x_{k+1}=x_k-\alpha_kg_k,
\end{equation}
where $g_k=\nabla f(x_k)$ and $\alpha_k>0$ is the stepsize. The classic choice for the stepsize is due to Cauchy \cite{cauchy1847methode} who suggested to calculate it by exact line search, i.e.,
\begin{equation}\label{alpsd0}
	\alpha_k^{SD}=\arg\min_{\alpha>0}~f(x_k-\alpha g_k), 
\end{equation}
which together with \eqref{eqitr} is often referred to as the steepest descent (SD) method. Many stepsizes have been designed to improve the efficiency of gradient methods such as the Barzilai--Borwein stepsize \cite{barzilai1988two}, the Polyak stepsize \cite{polyak1969minimization,polyak1987introduction} and the Yuan stepsize \cite{yuan2006new}. Due to the good performance in solving machine learning problems, there has been a surge of interest in the Polyak stepsize  \cite{jiang2024adaptive,loizou2021stochastic,wang2023generalized}, which  minimizes an upper bound of the distance between the new iteration to the optimal solution and has the form 
\begin{equation}\label{plks}
	\alpha_k^{P} = 2\gamma\frac{f_k-f_*}{\|g_k\|^2},\quad \gamma\in(0,1].
\end{equation}
where $f_k=f(x_k)$ and $f_*$ is the optimal function value. Convergence of the gradient method using the Polyak stepsize has been studied in different works for the case $\gamma=\frac{1}{2}$ or 1  \cite{barre2020complexity,hazan2019revisiting}.

In recent years, analyzing the worst-case complexities of gradient-based methods has triggered much attention \cite{cartis2010complexity,das2024branch,de2020worst,drori2014performance,gu2020tight,lessard2016analysis,taylor2017exact,teboulle2023elementary,yuan2010short}. Particularly, Drori and Teboulle \cite{drori2014performance} formulate the worst-case complexity of a given algorithm as an optimal solution to the so-called performance estimation problem (PEP).  By making some relaxations to PEP, a bound of the worst-case complexity can be obtained by solving a convex semidefinite programming (SDP). Following \cite{drori2014performance}, the worst-case complexities  for various popular methods have been analyzed including the gradient method, Nesterov's fast gradient method, inexact gradient method, proximal point algorithm, see  \cite{de2020worst,gu2020tight,taylor2017exact} and references therein.

Although the PEP framework is a powerful tool for the study of algorithm complexity, as pointed out by Teboulle and Vaisbourd \cite{teboulle2023elementary}, it lacks an intuitive explanation how the worst-case complexity is addressed: (i) one often has to employ a computer-assisted proof since the aforementioned SDP usually has no closed form solution and must be solved numerically; (ii) it is generally unclear how the solution is deduced when a closed form solution is available.  Moreover, it is not easy to apply the PEP framework to the gradient method with adaptive stepsizes \cite{goujaud2023fundamental}. Consequently, most of existing complexity results of gradient methods focus on fixed stepsizes. Recently, by resorting to PEP, the worst-case complexities of the gradient method with exact line search and the Polyak stepsize \eqref{plks} with $\gamma=1$ are presented through computer-assisted proof in \cite{de2017worst} and \cite{barre2020complexity}, respectively.


The purpose of this paper is to give a novel analytical analysis for the worst-case complexity of the gradient method with exact line search and the Polyak stepsize, respectively, for smooth strongly convex functions with Lipschitz continuous gradient. To this end, we first study a family of gradient methods whose stepsizes include the one determined by exact line search and the Polyak stepsize as special instances. We show that the family converges linearly when applying to a strongly convex quadratic function. Then, based on the convergence results on quadratics, for general smooth strongly convex objective functions, we establish the worst-case complexity of the gradient method with exact line search and the Polyak stepsize in an analytic way, respectively, which recovers the computer-assisted results in \cite{de2017worst} and \cite{barre2020complexity}. To our knowledge, this is the first analytic proof of the worst-case complexity of the gradient method with exact line search and the Polyak stepsize. Furthermore, we show that the considered family will asymptotically zigzag in a two-dimensional subspace spanned by the two eigenvectors corresponding to the largest and smallest eigenvalues of the Hessian, which has not been clarified in the literature for the gradient method using the Polyak stepsize.

The paper is organized as follows. In Section \ref{sec2}, we investigate linear convergence of a family of gradient methods on strongly convex quadratics. In Section \ref{sec3}, we present our analytic proof of the worst-case complexity of the gradient method with exact line search and the Polyak stepsize, respectively. Finally, some conclusion remarks are drawn and the asymptotic behavior of the considered family is discussed  in Section \ref{sec4}.

\section{Convergence analysis for quadratics}\label{sec2}
In this section, we consider the unconstrained quadratic optimization 
\begin{equation}\label{qudpro}
	\min_{x\in\mathbb{R}^n}~f(x)=\frac{1}{2}x\tr Ax-b\tr x,
\end{equation}
where $b\in \mathbb{R}^n$ and $A\in \mathbb{R}^{n\times n}$ is symmetric positive definite. Let $\{\lambda_1,\lambda_2,\cdots,\lambda_n\}$ be the eigenvalues of $A$, and $\{\xi_1,\xi_2,\ldots,\xi_n\}$ be the associated orthonormal eigenvectors.  Without loss of generality, we assume that
\begin{equation*}\label{formA}
	A=\textrm{diag}\{\lambda_1,\lambda_2,\cdots,\lambda_n\},~~0<\mu=\lambda_1<\lambda_2<\cdots<\lambda_n=L.
\end{equation*}

For a more general and uniform analysis, we investigate a family of gradient methods with stepsize given by
\begin{equation}\label{fyalp2}
	\alpha_k=\gamma\frac{g_{k}\tr   \psi(A) g_{k}}{g_{k}\tr   \psi(A)A g_{k}},\quad \gamma\in(0,1],
\end{equation}
where $\psi$ is a real analytic function on $[\mu,L]$ and can be expressed by Laurent series $\psi(z)=\sum_{k=-\infty}^\infty c_kz^k$ with $c_k\in\mathbb{R}$ such that $0<\sum_{k=-\infty}^\infty c_kz^k<+\infty$ for all $z\in[\mu,L]$. 
Apparently, when $\psi(A)=I$, we get the generalized SD stepsize, namely
\begin{equation}\label{alpsd}
	\alpha_k^{GSD}=\gamma\arg\min_{\alpha>0}~f(x_k-\alpha g_k)=\gamma\frac{g_{k} \tr  g_{k}}{g_{k} \tr  Ag_{k}},
\end{equation}
which is exactly the stepsize $\alpha_k^{SD}$ obtained by exact line search when  $\gamma=1$. See \cite{dai2003altermin} for gradient methods with shorten SD steps, i.e., $\gamma\in(0,1)$. Interestingly, the Polyak stepsize $\alpha_k^{P}$ in \eqref{plks} corresponding to the case $\psi(A)=A^{-1}$. In fact, since $f_*=f(x_*)$, where $x_*=A^{-1}b$ is the optimal solution, we have
\begin{equation*}
	f_k=f_*+\frac{1}{2}(x_k-x_*)\tr A(x_k-x_*)=f_*+\frac{1}{2}g_k\tr A^{-1}g_k,
\end{equation*}
which gives
\begin{equation}\label{plksqp}
	\alpha_k^{P}=2\gamma\frac{f_k-f_*}{\|g_k\|^2} =\gamma\frac{g_k\tr A^{-1}g_k}{g_k\tr g_k}.
\end{equation}

For a symmetric positive definite matrix $D$, denote the weighted norm by $\|x\|_{D}=\sqrt{x\tr Dx}$. We establish line convergence of the family of gradient methods \eqref{fyalp2} in the next theorem. 
\begin{theorem}\label{thrateqp}
	Consider applying the gradient method \eqref{eqitr} with the stepsize $\alpha_k$ given by \eqref{fyalp2} 	to solve the unconstrained quadratic optimization \eqref{qudpro}. It holds that
	\begin{align}\label{gdfrateqp}
		\|x_{k+1}-x_*\|_{\psi(A)A}^2
		&\leq \left(1-\frac{4\gamma(2-\gamma)\mu L}{(L+\mu)^2}\right)\|x_{k}-x_*\|_{\psi(A)A}^2.
	\end{align} 
	Moreover, if  $x_k=A^{-1}\left(\frac{1}{\sqrt{2}}\psi^{-\frac{1}{2}}(A)\left(\xi_1\pm\xi_n\right)+b\right)$,  then the two sides of \eqref{gdfrateqp} are equal.
\end{theorem} 
\begin{proof}
	Using $x_{k}-x_*=A^{-1}g_k$  and the definition of $\alpha_k$ in \eqref{fyalp2}, we have
	\begin{align}\label{gdfrateqp2}
		\|x_{k+1}-x_*\|_{\psi(A)A}^2&=	\|x_{k}-\alpha_kg_k-x_*\|_{\psi(A)A}^2\nonumber\\
		&=\|x_{k}-x_*\|_{\psi(A)A}^2-2\alpha_kg_k\tr \psi(A)A(x_{k}-x_*)+\alpha_k^2g_k\tr \psi(A)Ag_k\nonumber\\
		&=g_k\tr \psi(A)A^{-1}g_k-2\alpha_kg_k\tr \psi(A)g_k+\alpha_k^2g_k\tr \psi(A)Ag_k\nonumber\\
		&=g_k\tr \psi(A)A^{-1}g_k-\gamma(2-\gamma)\frac{(g_k\tr \psi(A)g_k)^2}{g_k\tr \psi(A)Ag_k}\nonumber\\
		&=\left(1-\frac{\gamma(2-\gamma)(g_k\tr \psi(A)g_k)^2}{g_k\tr \psi(A)Ag_k\cdot g_k\tr \psi(A)A^{-1}g_k}\right)g_k\tr \psi(A)A^{-1}g_k\nonumber\\
		&=\left(1-\frac{\gamma(2-\gamma)(g_k\tr \psi(A)g_k)^2}{g_k\tr \psi(A)Ag_k\cdot g_k\tr \psi(A)A^{-1}g_k}\right)\|x_{k}-x_*\|_{\psi(A)A}^2\nonumber\\
		&\leq \left(1-\frac{4\gamma(2-\gamma)\mu L}{(L+\mu)^2}\right)\|x_{k}-x_*\|_{\psi(A)A}^2,
	\end{align}
	where the last inequality follows from the Kantorovich inequality.

	If   $x_k=A^{-1}\left(\frac{1}{\sqrt{2}}\psi^{-\frac{1}{2}}(A)\left(\xi_1\pm\xi_n\right)+b\right)$, we have
	\begin{equation}\label{gkrateqp1}
		g_k=\frac{1}{\sqrt{2}}\psi^{-\frac{1}{2}}(A)\left(\xi_1\pm\xi_n\right),
	\end{equation}
	which gives
	\begin{equation*}
		g_k\tr \psi(A)A^{j}g_k=\frac{1}{2}\left(\mu^{j}+L^{j}\right),\quad j=-1,0,1.
	\end{equation*}
	Thus,
	\begin{align*}
		\|x_{k+1}-x_*\|_{\psi(A)A}^2
		&=\left(1-\frac{\gamma(2-\gamma)}{\frac{1}{2}(L+\mu)\cdot \frac{1}{2}\left(\mu^{-1}+L^{-1}\right)}\right)\|x_k-x_*\|_{\psi(A)A}^2\nonumber\\
		&= \left(1-\frac{4\gamma(2-\gamma)\mu L}{(L+\mu)^2}\right)\|x_k-x_*\|_{\psi(A)A}^2.
	\end{align*}
	We complete the proof.
\end{proof}

\begin{remark}\label{remk1}
	From the proof of Theorem \ref{thrateqp}, if the worst-case of convergence rate of the gradient method with the stepsize $\alpha_k$ given by \eqref{fyalp2} is achieved, we must have $\alpha_{k}=\frac{2\gamma}{L+\mu}$. Moreover, when $\gamma=1$ and $x_t=A^{-1}\left(\frac{1}{\sqrt{2}}\psi^{-\frac{1}{2}}(A)\left(\xi_1\pm\xi_n\right)+b\right)$ for some $t$,  then the two sides of \eqref{gdfrateqp} are equal for all $k\geq t$.
\end{remark}

%
%

By setting $\psi(A)=I$ and $A^{-1}$ in \eqref{gdfrateqp} respectively, we get the following results for the gradient method with the generalized SD stepsize  \eqref{alpsd} and the Polyak stepsize \eqref{plksqp}.
\begin{corollary}\label{corrateqp}
	Consider applying the gradient method \eqref{eqitr} to solve the unconstrained quadratic optimization \eqref{qudpro}. 
	
	(i) If $\alpha_k$ is given by the generalized SD stepsize  \eqref{alpsd}, then 
	\begin{equation}\label{sdrateqp}
		f_{k+1}-f_*
		\leq \left(1-\frac{4\gamma(2-\gamma)\mu L}{(L+\mu)^2}\right)(f_{k}-f_*).
	\end{equation} 
	
	(ii) If $\alpha_k$ is given by the Polyak stepsize \eqref{plksqp}, then 
	\begin{equation}\label{polyakrateqp}
		\|x_{k+1}-x_*\|^2
		\leq \left(1-\frac{4\gamma(2-\gamma)\mu L}{(L+\mu)^2}\right)\|x_{k}-x_*\|^2.
	\end{equation} 
		%
\end{corollary}

\begin{remark}\label{remk2}
	From Corollary \ref{corrateqp}, for the generalized SD stepsize \eqref{alpsd}, $\alpha_k^{SD}=\frac{g_{k} \tr  g_{k}}{g_{k} \tr  Ag_{k}}$ yields faster convergence rate than other values of $\gamma$ while for the Polyak stepsize \eqref{plksqp}, $\alpha_k^{P}=2\frac{f_k-f_*}{\|g_k\|^2}$ gives the fastest gradient method in the sense of \eqref{polyakrateqp}. 
\end{remark}

\section{Worst-case complexity for smooth strongly convex functions}\label{sec3}
In this section, by making use of convergence results in the former section, we present an analytic analysis of the worst-case complexity for the gradient method with exact line search and Polyak stepsize, respectively, on the class of $L$-smooth and $\mu$-strongly convex functions, denoted as $\mathcal{F}_{\mu,L}$.

\begin{definition}
	Let	$f: \mathbb{R}^n \to \mathbb{R}$ be a continuously differentiable function and $\mu,L>0$. We say that $f$ is  $L$-smooth and $\mu$-strongly convex if
	
	(i) $f$ is $L$-smooth, i.e.,
	\begin{equation*}
		f(y)\leq f(x)+\nabla f(x)\tr (y-x)+\frac{L}{2}\|x-y\|^2,~\forall~x,y\in\mathbb{R}^n,
	\end{equation*}	
	
	(ii) $f$ is $\mu$-strongly convex, i.e.,
	\begin{equation*}
		f(y)\geq f(x)+\nabla f(x)\tr (y-x)+\frac{\mu}{2}\|x-y\|^2,~\forall~x,y\in\mathbb{R}^n.
	\end{equation*}
\end{definition}

We will use the following important inequality for $f\in\mathcal{F}_{\mu,L}$, 
\begin{align}\label{intpscf}
	f(x)&-f(y)+\nabla f(x)\tr (y-x)+\frac{1}{2(1-\mu/L)}\times\nonumber\\
	&\left(\mu\|x-y\|^2-\frac{2\mu}{L}(\nabla f(x)-\nabla f(y))\tr (x-y)+\frac{1}{L}\|\nabla f(x)-\nabla f(y)\|^2\right)\leq0,
\end{align} 
see Theorem 4 of \cite{taylor2017smooth}.

\subsection{Gradient method with exact line search}

Now we show analytically the worst-case complexity of the gradient method with exact line search.

\begin{theorem}\label{thsdworst}
	Suppose that $f\in\mathcal{F}_{\mu,L}$ and consider the gradient method with stepsize $\alpha_k$ determined by exact line search. Then, 
	\begin{equation}\label{sdworstrate1}
		f_{k+1}-f_*\leq\left(\frac{L-\mu}{L+\mu}\right)^2(f_k-f_*).
	\end{equation}
\end{theorem}
\begin{proof}
	It follows from \eqref{intpscf}, $g_*=0$  and $g_{k+1}\tr g_{k}=0$ that
	\begin{align}
		&f_{k+1}-f_k+\frac{1}{2(1-\mu/L)}\nonumber\\
		&\times\left(\mu\|x_{k+1}-x_{k}\|^2+\frac{2\mu}{L}g_{k}\tr (x_{k+1}-x_{k})+\frac{1}{L}\|g_{k+1}-g_{k}\|^2\right)\leq0,\label{less0fkfk1}\\
		&f_k-f_*+g_{k}\tr (x_*-x_{k})+\frac{1}{2(1-\mu/L)}\nonumber\\
		&\quad\times\left(\mu\|x_{k}-x_*\|^2-\frac{2\mu}{L}g_{k}\tr (x_{k}-x_*)+\frac{1}{L}\|g_{k}\|^2\right)\leq0,\label{less0fkfs}
	\end{align} 
	and
	\begin{align}
		&f_{k+1}-f_*+g_{k+1}\tr (x_*-x_{k+1})+\frac{1}{2(1-\mu/L)}\nonumber\\
		&\times\left(\mu\|x_{k+1}-x_*\|^2-\frac{2\mu}{L}g_{k+1}\tr (x_{k+1}-x_*)+\frac{1}{L}\|g_{k+1}-g_*\|^2\right)\leq0.\label{less0fk1fs}
	\end{align} 
	Suppose that $\zeta_1,\zeta_2,\zeta_3$ are nonnegative scalars and $\zeta_1+\zeta_2+\zeta_3>0$. Then, weighted sum of the above three inequalities yields 
	\begin{align*}
		0&\geq\zeta_1\times\eqref{less0fkfk1}+\zeta_2\times\eqref{less0fkfs}+\zeta_3\times\eqref{less0fk1fs}\\
		&=(\zeta_1+\zeta_3)(f_{k+1}-f_*)+(\zeta_2-\zeta_1)(f_k-f_*)
		+\frac{1}{2(1-\mu/L)}\\
		&\times\Bigg[\mu\zeta_1\|(x_{k+1}-x_*)-(x_{k}-x_*)\|^2+\mu\zeta_2\|x_{k}-x_*\|^2+\mu\zeta_3\|x_{k+1}-x_*\|^2\\
		&+\frac{2\mu\zeta_1}{L}g_{k}\tr ((x_{k+1}-x_*)-(x_{k}-x_*))+2\zeta_2g_{k}\tr (x_*-x_{k})
		\\
		&+2\zeta_3g_{k+1}\tr (x_*-x_{k+1})
		+\frac{\zeta_1+\zeta_3}{L}\|g_{k+1}\|^2+\frac{\zeta_1+\zeta_2}{L}\|g_{k}\|^2\Bigg]\\
		&=(\zeta_1+\zeta_3)(f_{k+1}-f_*)+(\zeta_2-\zeta_1)(f_k-f_*)
		+\frac{1}{2(1-\mu/L)}\\
		&\times\Bigg[
		\mu(\zeta_1+\zeta_2)\|x_{k}-x_*\|^2-\left(\frac{2\mu\zeta_1}{L}+2\zeta_2\right)g_{k}\tr (x_{k}-x_*)\\
		&
		-2\mu\zeta_1(x_{k+1}-x_*)\tr (x_{k}-x_*)+\left(\frac{2\mu\zeta_1}{L}g_{k}-2\zeta_3g_{k+1}\right)\tr (x_{k+1}-x_*)
		\\
		&
		+\mu(\zeta_1+\zeta_3)\|x_{k+1}-x_*\|^2+\frac{\zeta_1+\zeta_3}{L}\|g_{k+1}\|^2+\frac{\zeta_1+\zeta_2}{L}\|g_{k}\|^2\Bigg].
	\end{align*}
	To proceed, we complete the square of all items containing $x_{k}-x_*$ to get 
	\begin{align*}
		0&\geq
		(\zeta_1+\zeta_3)(f_{k+1}-f_*)+(\zeta_2-\zeta_1)(f_k-f_*)
		+\frac{1}{2(1-\mu/L)}\\
		&\times\Bigg[
		\mu(\zeta_1+\zeta_2)
		\left\|x_{k}-x_*-\frac{\zeta_1}{\zeta_1+\zeta_2}(x_{k+1}-x_*)-\frac{\mu\zeta_1+L\zeta_2}{\mu L(\zeta_1+\zeta_2)}g_{k}-\delta g_{k+1}\right\|^2\\
		&+\beta_1\|x_{k+1}-x_*\|^2+\beta_2\|g_{k}\|^2
		+\beta_3\|g_{k+1}\|^2
		-2(x_{k+1}-x_*)\tr \left(\beta_4g_{k}-\beta_5g_{k+1}\right)
		\Bigg]
		\\
		&=(\zeta_1+\zeta_3)(f_{k+1}-f_*)+(\zeta_2-\zeta_1)(f_k-f_*)
		+\frac{1}{2(1-\mu/L)}\\
		&\times\Bigg[
		\mu(\zeta_1+\zeta_2)
		\left\|x_{k}-x_*-\frac{\zeta_1}{\zeta_1+\zeta_2}(x_{k+1}-x_*)-\frac{\mu\zeta_1+L\zeta_2}{\mu L(\zeta_1+\zeta_2)}g_{k}-\delta g_{k+1}\right\|^2\\
		&+\beta_1\left\|x_{k+1}-x_*-\frac{1}{\beta_1}\left(\beta_4g_{k}-\beta_5g_{k+1}
		\right)\right\|^2\\
		&+\left(\beta_2-\frac{\beta_4^2}{\beta_1}\right)\|g_{k}\|^2+\left(\beta_3-\frac{\beta_5^2}{\beta_1}\right)\|g_{k+1}\|^2\Bigg],
	\end{align*} 
	where $\delta$ is some scalar and
	\begin{equation*}
		\beta_1=\mu(\zeta_1+\zeta_3)-\frac{\mu\zeta_1^2}{\zeta_1+\zeta_2},\quad
		\beta_2=\frac{(L-\mu)(\mu\zeta_1^2-L\zeta_2^2)}{\mu L^2(\zeta_1+\zeta_2)}, 
	\end{equation*}
	\begin{equation*}
		\beta_3=\frac{\zeta_1+\zeta_3}{L}-\mu(\zeta_1+\zeta_2) \delta^2,\quad
		\beta_4=\frac{(L-\mu)\zeta_1\zeta_2}{L(\zeta_1+\zeta_2)},\quad 
		\beta_5=\mu\delta \zeta_2-\zeta_3.
	\end{equation*}
	As we will see later $\delta\neq0$ is important to derive the worst-case rate. It follows that
	\begin{align}
		&f_{k+1}-f_*\leq
		\frac{\zeta_2-\zeta_1}{\zeta_1+\zeta_3}(f_k-f_*)
		-\frac{1}{2(1-\mu/L)}\nonumber\\
		&\times\Bigg[
		\frac{\mu(\zeta_1+\zeta_2)}{\zeta_1+\zeta_3}
		\left\|x_{k}-x_*-\frac{\zeta_1}{\zeta_1+\zeta_2}(x_{k+1}-x_*)-\frac{\mu\zeta_1+L\zeta_2}{\mu L(\zeta_1+\zeta_2)}g_{k}-\delta g_{k+1}\right\|^2\nonumber\\
		&+\frac{\beta_1}{\zeta_1+\zeta_3}\left\|x_{k+1}-x_*-\frac{1}{\beta_1}\left(\beta_4g_{k}-\beta_5g_{k+1}
		\right)\right\|^2\nonumber\\
		&+\frac{1}{\zeta_1+\zeta_3}\left(\beta_2-\frac{\beta_4^2}{\beta_1}\right)\|g_{k}\|^2+\frac{1}{\zeta_1+\zeta_3}\left(\beta_3-\frac{\beta_5^2}{\beta_1}\right)\|g_{k+1}\|^2\Bigg].\label{sdconreq1}
	\end{align}

	With the aim of determining $\zeta_1$, $\zeta_2$ and $\zeta_3$, we consider applying the gradient method with exact line search to a two-dimensional quadratic function with 
		\begin{equation}\label{sd2dwsr}
			A=\begin{pmatrix}
				\mu & 0 \\
				0  & L \\
			\end{pmatrix}.
		\end{equation}
		By Theorem \ref{thrateqp} and \eqref{sdrateqp}, in this case, the worst-case rate of the gradient method with exact line search  matches the one in \eqref{sdworstrate1}.  
		Since the inequality \eqref{sdconreq1} applies to the above quadratic case and the worst-case rate in \eqref{sdrateqp} does not involve the last four terms in \eqref{sdconreq1}, we require
		\begin{equation}\label{betas0eq1}
			\beta_2-\frac{\beta_4^2}{\beta_1} =0,\quad \beta_3-\frac{\beta_5^2}{\beta_1}=0,
		\end{equation}
		\begin{equation}\label{xkgksquare1}
			\left\|x_{k}-x_*-\frac{\zeta_1}{\zeta_1+\zeta_2}(x_{k+1}-x_*)-\frac{\mu\zeta_1+L\zeta_2}{\mu L(\zeta_1+\zeta_2)}g_{k}-\delta g_{k+1}\right\|^2=0,
		\end{equation}
		and
		\begin{equation}\label{xkgksquare2}
			\left\|x_{k+1}-x_*-\frac{1}{\beta_1}\left(\beta_4g_{k}-\beta_5g_{k+1}
			\right)\right\|^2=0,
		\end{equation}
		From Remark \ref{remk1}, we must have $\alpha_{k}=\frac{2}{L+\mu}$ when achieving the worst-case rate. Thus, for the above two-dimensional quadratic function, we get
		\begin{equation*}
			x_{k}-x_*=A^{-1}g_k,\quad
			g_{k+1}=(I-\alpha_{k}A)g_k=\begin{pmatrix}
				1-\frac{2\mu}{L+\mu} & 0 \\
				0  & 1-\frac{2L}{L+\mu} \\
			\end{pmatrix}g_k.
		\end{equation*}
		Suppose that the two components of $g_k$ are nonzero. The above relations together with \eqref{xkgksquare1} and \eqref{xkgksquare2} imply that 
		\begin{equation}\label{xkgk1s0qp1}
			\frac{1}{\mu}-\left(\frac{\zeta_1}{\zeta_1+\zeta_2}\frac{1}{\mu}+\delta\right)\left(1-\frac{2\mu}{L+\mu}\right)-\frac{\mu\zeta_1+L\zeta_2}{\mu L(\zeta_1+\zeta_2)}=0,
		\end{equation}
		\begin{equation}\label{xkgk1s0qp2}
			\frac{1}{L}-\left(\frac{\zeta_1}{\zeta_1+\zeta_2}\frac{1}{L}+\delta\right)\left(1-\frac{2L}{L+\mu}\right)-\frac{\mu\zeta_1+L\zeta_2}{\mu L(\zeta_1+\zeta_2)}=0,
		\end{equation}
		\begin{equation}\label{betas0qp1}
			\left(\frac{\beta_1}{\mu}+\beta_5\right)\left(1-\frac{2\mu}{L+\mu}\right)-\beta_4=0,
		\end{equation}
		\begin{equation}\label{betas0qp2}
			\left(\frac{\beta_1}{L}+\beta_5\right)\left(1-\frac{2L}{L+\mu}\right)-\beta_4=0.
		\end{equation}
		Then, $\mu L\times(\eqref{xkgk1s0qp1}+\eqref{xkgk1s0qp2})$ yields
		\begin{equation*}
			L+\mu-\frac{\zeta_1}{\zeta_1+\zeta_2}\frac{(L-\mu)^2}{L+\mu}-2\frac{\mu\zeta_1+L\zeta_2}{(\zeta_1+\zeta_2)}=0,
		\end{equation*}
		which indicates
		\begin{equation}\label{zt1zt2rel1}
			\zeta_2=\frac{2\mu}{L+\mu}\zeta_1.
		\end{equation}
		From the first equation in \eqref{betas0eq1}, we get
		\begin{equation*}
			\beta_1=\frac{\beta_4^2}{\beta_2} =\frac{(L-\mu)^2\zeta_1^2\zeta_2^2}{L^2(\zeta_1+\zeta_2)^2}\frac{\mu L^2(\zeta_1+\zeta_2)}{(L-\mu)(\mu\zeta_1^2-L\zeta_2^2)}=\frac{\mu(L-\mu)\zeta_1^2\zeta_2^2}{(\zeta_1+\zeta_2)(\mu\zeta_1^2-L\zeta_2^2)}.
		\end{equation*}
		Using the definition of $\beta_1$ and rearranging terms, we obtain
		\begin{equation}\label{zt1zt3rel1}
			\zeta_3=\frac{\mu\zeta_1^2(\zeta_1-\zeta_2)}{\mu\zeta_1^2-L\zeta_2^2}-\zeta_1=\frac{2\mu}{L-\mu}\zeta_1.
		\end{equation}
		By \eqref{xkgk1s0qp1} or \eqref{xkgk1s0qp2}, we have
		\begin{equation}\label{zt1dltrel1}
			\delta=\frac{1}{L(\zeta_1+\zeta_2)}\zeta_1=\frac{L+\mu}{L(L+3\mu)}.
		\end{equation}
		It is not difficult to check that the choices $\zeta_2$, $\zeta_3$ and $\delta$ in \eqref{zt1zt2rel1}, \eqref{zt1zt3rel1} and \eqref{zt1dltrel1} are such that the two equations in \eqref{betas0eq1}, \eqref{betas0qp1} and \eqref{betas0qp2} hold. Thus, we only need to find nonnegative $\zeta_1$, $\zeta_2$ and $\zeta_3$ satisfying \eqref{zt1zt2rel1} and \eqref{zt1zt3rel1}. Letting $\zeta_3=1$, we get
		\begin{equation}\label{parassd1}
			\zeta_1=\frac{L-\mu}{2\mu},\quad \zeta_2=\frac{L-\mu}{L+\mu},\quad \zeta_3=1.
		\end{equation}

		Substituting the above $\zeta_1,\zeta_2,\zeta_3$ in \eqref{parassd1} into \eqref{sdconreq1}, we get
		\begin{align*}
			&f_{k+1}-f_*
			\leq\left(\frac{L-\mu}{L+\mu}\right)^2(f_k-f_*)
			\\
			&
			-\frac{\mu L(L+3\mu)}{2(L+\mu)^2}
			\left\|x_{k}-x_*-\frac{L+\mu}{L+3\mu}(x_{k+1}-x_*)-\frac{3L+\mu}{L(L+3\mu)}g_{k}-\frac{L+\mu}{L(L+3\mu)}g_{k+1}\right\|^2\\
			&-\frac{2\mu L^2}{(L-\mu)(L+3\mu)}\left\|x_{k+1}-x_*-\frac{(L-\mu)^2}{2\mu L(L+\mu)}g_{k}-\frac{L+\mu}{2\mu L}g_{k+1}
			\right\|^2,
		\end{align*}	
		which implies \eqref{sdworstrate1}. In addition, the rate in \eqref{sdworstrate1} is achieved by applying the gradient method with exact line search to a two-dimensional quadratic function with Hessian given by \eqref{sd2dwsr}. We complete the proof.
	\end{proof}

	\begin{remark}
		If we let $\zeta_1+\zeta_3=1$, it follows from \eqref{zt1zt2rel1} and \eqref{zt1zt3rel1} that
		\begin{equation*}
			\zeta_1=\frac{L-\mu}{L+\mu},\quad \zeta_2=\frac{2\mu(L-\mu)}{(L+\mu)^2},\quad \zeta_3=\frac{2\mu}{L+\mu},
		\end{equation*}
		which are the parameters in \cite{de2017worst}.
	\end{remark}

	\subsection{Gradient method with the Polyak stepsize}
	
	Next theorem gives an analytic analysis of the worst-case complexity of the gradient method using the Polyak stepsize \eqref{plks}. 
	
	\begin{theorem}\label{thpolyak}
		Suppose that $f\in\mathcal{F}_{\mu,L}$ and consider the gradient method with stepsize $\alpha_k$ given by \eqref{plks}. 	Then, 
		\begin{equation}\label{wstratePolyak}
			\|x_{k+1}-x_*\|^2\leq\left(1-\frac{4\gamma(2-\gamma)\mu L}{(L+\mu)^2}\right)\|x_{k}-x_*\|^2.
		\end{equation}
	\end{theorem}
	\begin{proof}
		By \eqref{intpscf} and $g_*=0$, we have 
		\begin{align}\label{fk2fsscl}
			f_k&-f_*+g_k\tr (x_*-x_{k})+\frac{1}{2(1-\mu/L)}\times\nonumber\\
			&\left(\mu\|x_{k}-x_*\|^2-\frac{2\mu}{L}g_k\tr (x_{k}-x_*)+\frac{1}{L}\|g_k\|^2\right)\leq0
		\end{align} 
		and
		\begin{align}\label{fs2fkscl}
			f_*&-f_k+\frac{1}{2(1-\mu/L)}
			\left(\mu\|x_{k}-x_*\|^2+\frac{2\mu}{L}g_k\tr (x_*-x_{k})+\frac{1}{L}\|g_k\|^2\right)\leq0.
		\end{align} 
		It follows from the definition of $\alpha_k$ in \eqref{plks} that
		\begin{align}\label{fs2fkeqalp}
			&2\gamma(f_k-f_*)-\alpha_k\|g_k\|^2=0.
		\end{align} 
		Let  $\zeta_1, \zeta_2$ be two nonnegative scalars, $\zeta_3\in \mathbb{R}$ and $\zeta_1+\zeta_2>0$. The following weighted sum is valid:
		\begin{align*}
			0&\geq\zeta_1\times\eqref{fk2fsscl}+\zeta_2\times\eqref{fs2fkscl}+\zeta_3\times\eqref{fs2fkeqalp}\\
			&=\zeta_1\Bigg[f_k-f_*+g_k\tr (x_*-x_{k})+\frac{1}{2(1-\mu/L)}\\
			&\times\left(\mu\|x_{k}-x_*\|^2-\frac{2\mu}{L}g_k\tr (x_{k}-x_*)+\frac{1}{L}\|g_k\|^2\right)\Bigg]\\
			&+\zeta_2\Bigg[f_*-f_k+\frac{1}{2(1-\mu/L)}
			\left(\mu\|x_{k}-x_*\|^2+\frac{2\mu}{L}g_k\tr (x_*-x_{k})+\frac{1}{L}\|g_k\|^2\right)\Bigg]\\
			&+\zeta_3\left[2\gamma(f_k-f_*)-\alpha_k\|g_k\|^2\right]\\
			&=(\zeta_1-\zeta_2+2\gamma\zeta_3)(f_k-f_*)+\mu L\beta\|x_{k}-x_*\|^2\\
			&-\left(2\mu\beta+\zeta_1\right)g_k\tr (x_{k}-x_*)+\left(\beta-\zeta_3\alpha_k\right)\|g_k\|^2,
		\end{align*}
		where $\beta=\frac{\zeta_1+\zeta_2}{2(L-\mu)}$. 	Substituting 
		\begin{equation*}
			\|x_{k+1}-x_*\|^2-\|x_{k}-x_*\|^2-\alpha_k^2\|g_k\|^2=-2\alpha_kg_k\tr (x_{k}-x_*)
		\end{equation*}
		into the above inequality, we obtain
		\begin{align}\label{rateeq1}
			&\sigma(f_k-f_*)+\frac{2\mu\beta+\zeta_1}{2\alpha_k}\|x_{k+1}-x_*\|^2\nonumber\\
			&+\left(\beta-\zeta_3\alpha_k-\frac{\alpha_k(2\mu\beta+\zeta_1)}{2}\right)\|g_k\|^2\leq\left(\frac{2\mu\beta+\zeta_1}{2\alpha_k}-\mu L\beta\right)\|x_{k}-x_*\|^2,
		\end{align}
		where $\sigma=\zeta_1-\zeta_2+2\gamma\zeta_3$.

		To get the desired convergence result, we assume that 
		\begin{equation}\label{gdnqeq0}
			\beta-\zeta_3\alpha_k-\frac{\alpha_k(2\mu\beta+\zeta_1)}{2}=0,
		\end{equation}
		which together with the definition of $\beta$ gives 
		\begin{equation*}
			\zeta_3=\frac{(1-\mu\alpha_k)(\zeta_1+\zeta_2)}{2(L-\mu)\alpha_k}-\frac{\zeta_1}{2}
		\end{equation*}
		and hence
		\begin{align}\label{zeta123t0}
			\sigma
			&=\zeta_1-\zeta_2+\frac{\gamma(1-\mu\alpha_k)}{(L-\mu)\alpha_k}(\zeta_1+\zeta_2)-\gamma\zeta_1
			\nonumber\\
			&=\frac{(1-\gamma)(L-\mu)\alpha_k +\gamma(1-\mu\alpha_k)}{(L-\mu)\alpha_k}\zeta_1+\frac{\gamma(1-\mu\alpha_k)-(L-\mu)\alpha_k}{(L-\mu)\alpha_k}\zeta_2.
		\end{align} 	
		Since $f\in\mathcal{F}_{\mu,L}$, we have
		\begin{equation*}\label{plksbd}
			\frac{\gamma}{L}\leq\alpha_k\leq\frac{\gamma}{\mu}.
		\end{equation*}
		Thus, $1-\mu\alpha_k\geq0$. Now we consider two cases:
		
		\textbf{Case 1}: $\gamma(1-\mu\alpha_k)-(L-\mu)\alpha_k<0$. In this case, we require
		\begin{equation}\label{fvleq0}
			\sigma=\zeta_1-\zeta_2+2\gamma\zeta_3=0,
		\end{equation}
		which together with \eqref{zeta123t0} indicates
		\begin{equation*}\label{zeta2t1}
			\zeta_2=\frac{(1-\gamma)(L-\mu)\alpha_k +\gamma(1-\mu\alpha_k)}{(L-\mu)\alpha_k-\gamma(1-\mu\alpha_k)}\zeta_1
		\end{equation*}
		and hence
		\begin{equation*}\label{betaval1}
			\beta=\frac{(2-\gamma)\alpha_k}{2((L-\mu)\alpha_k-\gamma(1-\mu\alpha_k))}\zeta_1.
		\end{equation*}
		Let
		\begin{equation}\label{zeta1val1}
			\zeta_1=2\frac{(L-\mu)\alpha_k-\gamma(1-\mu\alpha_k)}{(2-\gamma)\alpha_k}.
		\end{equation}
		We obtain $\beta=1$ and
		\begin{equation}\label{zeta3betaval1}
			\zeta_2=2\frac{(1-\gamma)(L-\mu)\alpha_k +\gamma(1-\mu\alpha_k)}{(2-\gamma)\alpha_k},\quad
			\zeta_3
			=\frac{2-(L+\mu)\alpha_k}{(2-\gamma)\alpha_k}.
		\end{equation}
		Clearly, $\zeta_1, \zeta_2$ are nonnegative. By \eqref{rateeq1}, \eqref{gdnqeq0}, \eqref{fvleq0}, \eqref{zeta1val1}   and \eqref{zeta3betaval1}, we get
		\begin{align}\label{rateeq2}
			\|x_{k+1}-x_*\|^2
			&\leq\left(1-\frac{2\mu L\beta\alpha_k}{2\mu\beta+\zeta_1}\right)\|x_{k}-x_*\|^2\nonumber\\
			&\leq\left(1-\frac{\mu L(2-\gamma)\alpha_k^2}{(L+\mu)\alpha_k-\gamma}\right)\|x_{k}-x_*\|^2.
		\end{align}
		Notice that the function
		$h(\alpha)=\frac{\alpha^2}{(L+\mu)\alpha-\gamma}$ achieves its minimum at $\frac{2\gamma}{L+\mu}$. 
		It follows from \eqref{rateeq2} that
		\begin{align*}
			\|x_{k+1}-x_*\|^2
			&\leq\left(1-(2-\gamma)\mu L h\left(\frac{2\gamma}{L+\mu}\right)\right)\|x_{k}-x_*\|^2\nonumber\\
			&=\left(1-\frac{4\gamma(2-\gamma)\mu L}{(L+\mu)^2}\right)\|x_{k}-x_*\|^2.
		\end{align*}

		\textbf{Case 2}: $\gamma(1-\mu\alpha_k)-(L-\mu)\alpha_k\geq0$. Recall that $\zeta_1, \zeta_2$ are nonnegative and $\zeta_1+\zeta_2>0$. In this case, we cannot expect \eqref{fvleq0} holds. 	It follows from \eqref{rateeq1}, \eqref{gdnqeq0} and \eqref{zeta123t0} that
		\begin{align}\label{rateeq4}
			\|x_{k+1}-x_*\|^2
			&\leq\left(1-\frac{2\mu L\beta\alpha_k}{2\mu\beta+\zeta_1}\right)\|x_{k}-x_*\|^2-\sigma(f_k-f_*).
		\end{align}
		To simplify our analysis, we set $\zeta_2=0$. By the strongly convexity of $f$,  the definition of $\beta$, \eqref{zeta123t0} and \eqref{rateeq4}, we have
		\begin{align*}
			\|x_{k+1}-x_*\|^2
			&\leq\left(1-\mu\alpha_k-\frac{\sigma\mu}{2}\right)\|x_{k}-x_*\|^2\nonumber\\
			&=\left(1-\mu\alpha_k-\mu\frac{(1-\gamma)(L-\mu)\alpha_k +\gamma(1-\mu\alpha_k)}{2(L-\mu)\alpha_k}\zeta_1\right)\|x_{k}-x_*\|^2.
		\end{align*}
		We further suppose that $\zeta_1$ is a scalar independent of $\alpha_k$. The function 
		\begin{equation*}
			h(\alpha)=\alpha+\frac{(1-\gamma)(L-\mu)\alpha+\gamma(1-\mu\alpha)}{2(L-\mu)\alpha}\zeta_1
		\end{equation*}
		attains its minimum at $\alpha=\sqrt{\frac{\gamma\zeta_1}{2(L-\mu)}}$, which implies 
		\begin{align}\label{rateeq6}
			\|x_{k+1}-x_*\|^2
			&\leq h\left(\sqrt{\frac{\gamma\zeta_1}{2(L-\mu)}}\right)\|x_{k}-x_*\|^2.
		\end{align}
		Notice that the rate in \eqref{rateeq6} holds for any $f\in\mathcal{F}_{\mu,L}$. So, for the quadratic case, the rate in \eqref{rateeq6} should coincide with the one in \eqref{polyakrateqp}  and the corresponding stepsize must be equal to $\frac{2\gamma}{L+\mu}$. This gives $\zeta_1=\frac{8\gamma(L-\mu)}{(L+\mu)^2}$. 	We get  \eqref{wstratePolyak} by substituting this $\zeta_1$ into \eqref{rateeq6}.

		We complete the proof by considering the above two cases.
	\end{proof}

	\begin{remark}
		When $\gamma=1$, the inequality \eqref{wstratePolyak} reduces to
		\begin{equation*}
			\|x_{k+1}-x_*\|^2\leq\left(\frac{L-\mu}{L+\mu}\right)^2\|x_{k}-x_*\|^2,
		\end{equation*}
		which recovers the one in Proposition 1 of \cite{barre2020complexity}. 
	\end{remark}

	\section{Concluding remarks}\label{sec4}
	

	Based on the convergence results of the family of gradient methods \eqref{fyalp2} on strongly convex quadratics, we have established the worst-case complexity of the gradient method with exact line search and the Polyak stepsize, respectively, in an analytic way for general smooth strongly convex objective functions. It is shown by Corollary \ref{corrateqp}  and Theorem \ref{thpolyak} that, from the worst-case complexity point of view, the gradient method using the Polyak stepsize \eqref{plksqp} achieves the fastest convergence rate when $\gamma=1$,  which is also true for the generalized SD stepsize \eqref{alpsd}. However, for the case $\gamma=1$, we can show that the family of gradient methods \eqref{fyalp2} will zigzag in a two-dimensional subspace spanned by the two eigenvectors corresponding to the largest and smallest eigenvalues of the Hessian. 
	Denoting the components of $g_k$ along the eigenvectors $\xi_i$ by $\nu_k^{(i)}$, $i=1,\ldots,n$, i.e., $g_k=\sum_{i=1}^n\nu_k^{(i)}\xi_i$. Next theorem gives the asymptotic convergence behavior of each gradient method in the family \eqref{fyalp2}, see Theorem 1 of \cite{huang2022asymptotic} for details.
	\begin{theorem}\label{thasconverg1}
		Assume that the starting point $x_0$ has the property that
		\begin{equation*}\label{assp1}
			g_0 \tr  \xi_1\neq0~~\textrm{and}~~g_0 \tr  \xi_n\neq0.
		\end{equation*}
		Let $\{x_k\}$ be the iterations generated by applying a method in \eqref{fyalp2} with $\gamma=1$ to solve  the unconstrained quadratic optimization \eqref{qudpro}. Then
		\begin{equation*}\label{mu2k}
			\lim_{k\rightarrow\infty}\frac{(\nu_{2k}^{(i)})^2}{\sum_{j=1}^n(\nu_{2k}^{(j)})^2}=
			\left\{
			\begin{array}{ll}
				\displaystyle\frac{1}{1+c^2}, & \hbox{if $i=1$,} \\
				0, & \hbox{if $i=2,\ldots,n-1$,} \\
				\displaystyle\frac{c^2}{1+c^2}, & \hbox{if $i=n$,}
			\end{array}
			\right.
		\end{equation*}
		and
		\begin{equation*}\label{mu2k1}
			\lim_{k\rightarrow\infty}\frac{(\nu_{2k+1}^{(i)})^2}{\sum_{j=1}^n(\nu_{2k+1}^{(j)})^2}=
			\left\{
			\begin{array}{ll}
				\displaystyle\frac{c^2\psi^2(L)}{\psi^2(\mu)+c^2\psi^2(L)}, & \hbox{if $i=1$,} \\
				0, & \hbox{if $i=2,\ldots,n-1$,} \\
				\displaystyle\frac{\psi^2(\mu)}{\psi^2(\mu)+c^2\psi^2(L)}, & \hbox{if $i=n$,}
			\end{array}
			\right.
		\end{equation*}
		where $c$ is a nonzero constant and can be determined by 
		\begin{equation*}\label{eqcvalue}
			c=\lim_{k\rightarrow\infty}\frac{\nu_{2k}^{(n)}}{\nu_{2k}^{(1)}}
			=-\frac{\psi(\mu)}{\psi(L)}\lim_{k\rightarrow\infty}\frac{\nu_{2k+1}^{(1)}}{\nu_{2k+1}^{(n)}}.
		\end{equation*}
	\end{theorem}	
	
	\begin{figure}[ht!b]
		\centering
		\includegraphics[width=0.75\textwidth,height=0.42\textwidth]{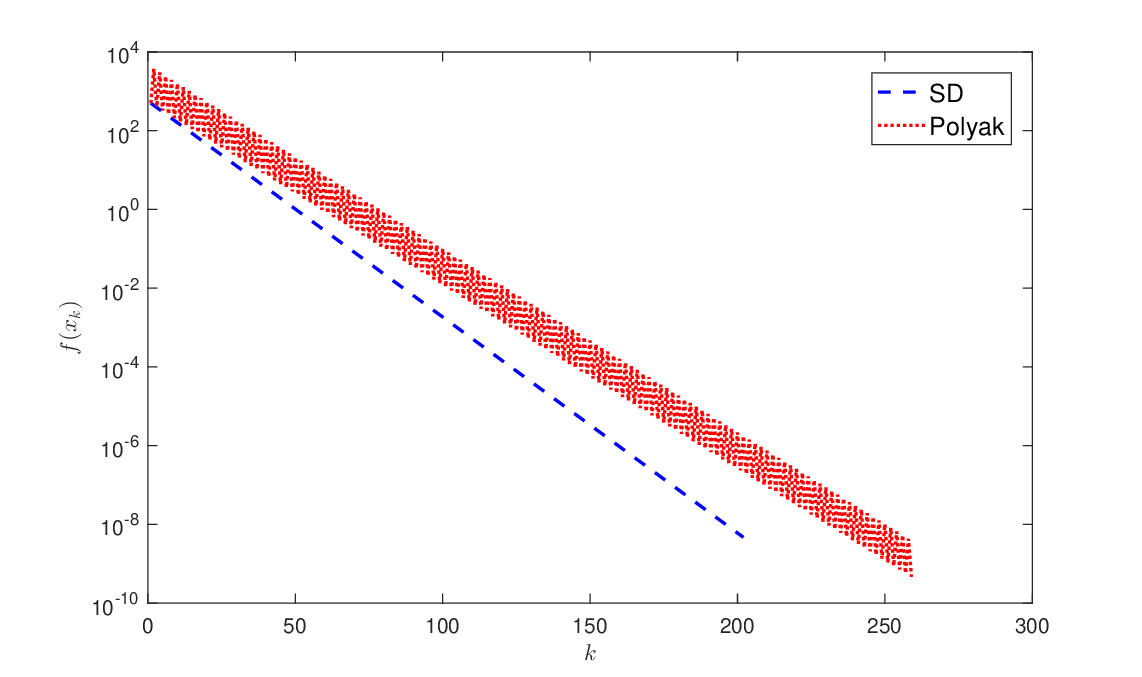}\\
		\caption{The gradient method with exact line search (SD) vs. that with the Polyak stepsize (Polyak).}\label{sdvspolyak}
	\end{figure}
	Due to the zigzag behavior shown in Theorem \ref{thasconverg1}, the gradient method using the Polyak stepsize \eqref{plksqp} with $\gamma=1$ may perform as poor as that with the exact line search. To see this, we apply the family of gradient methods \eqref{fyalp2} with $\gamma=1$ for $\psi(A)=I$ (exact line search) and $A^{-1}$ (the Polyak stepsize), respectively,  to minimize the quadratic problem \eqref{qudpro} with
	\begin{equation*}\label{tp1}
		A=\textrm{diag}\{1,100\} \quad \mbox{and} \quad b=0.
	\end{equation*}
	We use $x_0=(30,1)\tr $ as the starting point and stop the iteration if $\|g_k\|\leq10^{-8}\|g_0\|$. 
	It can be seen from Figure \ref{sdvspolyak} that the gradient method  with the Polyak stepsize may be slower than the SD method even for a strongly convex quadratic function. Techniques to break the zigzag pattern of the family of gradient methods \eqref{fyalp2} with $\gamma=1$ have been developed in \cite{huang2022asymptotic} which show promising performance on the gradient method with exact line search. We are wondering whether the techniques in \cite{huang2022asymptotic} can lead to efficient gradient methods based on the Polyak stepsize.

\begin{acknowledgements}
The work of the first author was supported by  Natural Science Foundation of Hebei Province (A2021202010). The work of the second author was supported by  Hong Kong RGC General Research Fund (15309223) and PolyU AMA Project (230413007).
\end{acknowledgements}


\bibliographystyle{siamplain}

\end{document}